\documentclass[10pt]{article}
\usepackage{times}
\usepackage{amsthm}
\usepackage{amssymb}
\usepackage{latexsym}
\usepackage{amscd}
\usepackage{epsfig}
\usepackage{verbatim}
\newtheorem{thm}{Theorem}
\newtheorem{que}{Question}[section]

\newcommand{\bd}{\ensuremath{\partial}}
\newcommand{\irr}{ir\-re\-duc\-i\-ble}

\newcommand{\tm}{3-man\-i\-fold}

\newcommand{\hm}{homeomorphic}
\newcommand{\Inte}{Int \, }

\newcommand{\p}{^{\prime}}

\begin{document}

\title{\normalsize \textbf{ON EMBEDDING THE FUNDAMENTAL GROUP OF 
A 3-MANIFOLD \\ IN ONE OF ITS KNOT GROUPS}}
\author{\small ROBERT MYERS \\ \small \textit{Department of Mathematics} 
\\ \small \textit{Oklahoma State University} \\  
\small \textit{Stillwater, OK 74078 USA} \\ 
\small \textit{Email: myersr@math.okstate.edu}}
\date{}

\maketitle
\begin{abstract} 
This paper gives necessary and sufficient conditions 
on a compact, connected, orientable 3-manifold $M$ 
for it to contain a knot $K$ such that $M-K$ is irreducible 
and $\pi_1(M)$ embeds in $\pi_1(M-K)$. 
This result provides counterexamples to a conjecture 
of Lopes and Morales and characterizes those 
orientable 3-manifolds for which it is true. 
\vspace{\baselineskip}

\noindent \textit{Keywords:} 3-manifold, fundamental group, knot group, 
irreducible, $\mathbb{Q}$ conjecture.

\medskip

\noindent Mathematics Subject Classification 2000: 57N10, 57M05 

\end{abstract}

\bigskip

\section*{\normalsize Introduction}
\stepcounter{section}

Given any connected \tm\ $M$ there is an easy way to construct a knot $K$ in 
$M$ such that $\pi_1(M)$ embeds in $\pi_1(M-K)$. Choose any 3-ball 
$B$ in $\Inte M=M-\bd M$, and choose any knot $K$ in $\Inte B$. 
Then $\pi_1(M-K)$ is isomorphic to $\pi_1(M)*\pi_1(B-K)$. 
Note that $M-K$ is reducible unless $M$ is \hm\ to $S^3$. 

Recently Lopes and Morales \cite{LM} have considered the question 
of finding a $K$ such that $M-K$ is \irr. They construct an infinite 
collection of compact, connected, orientable \tm s $M$ for which  
there is a two-component link $L$ in $M$ such that $M-L$ is \irr\ and 
$\pi_1(M)$ embeds in $\pi_1(M-L)$. A closer examination of these 
examples shows that each of these $M$ contains a knot $K$ such  
that $M-K$ is irreducible and  
$\pi_1(M)$ embeds in $\pi_1(M-K)$. (See section 3 below.)  

They conjecture \cite[Conjecture 1.3]{LM} that every compact \tm\ contains 
such a knot. 
If, as is customary and is the case in the examples above, one requires 
$K$ to lie in $\Inte M$, then any $M$ whose boundary contains 
a 2-sphere is a counterexample, since $\Inte M$ then contains no 
knots with \irr\ complement. Also, any $M$ such that $\pi_1(M)$ 
has torsion is a counterexample since the exterior 
of $K$ would be a Haken manifold and so $\pi_1(M-K)$ would contain 
no torsion. These observations led the author to 
wonder for which $M$ their conjecture is true. The following 
result gives an answer in the orientable case.   

\begin{thm} Let $M$ be a compact, connected, orientable 3-manifold. 
The following are equivalent.  
\begin{enumerate} 
\item If $M_i$ is a prime connected summand of $M$, 
then $M_i$ either 
\begin{enumerate} 
\item is homeomorphic to $S^1 \times S^2$, 
\item is a homotopy 3-sphere, or 
\item has non-empty boundary and is not a 3-ball. 
\end{enumerate}
\item $\bd M$ contains no 2-spheres; $\Inte M$ contains a 
compact, connected, irreducible 3-manifold $V$ with $\bd V\neq\emptyset$  
such that $\pi_1(V)\rightarrow\pi_1(M)$ is an isomorphism, 
and at most one component of $M-\Inte V$ is not 
a collar on a component of $\bd M$. 
\item There is a knot $K$ in $\Inte M$ such that $M-K$ 
is irreducible and $\pi_1(M)$ embeds in $\pi_1(M-K)$
\end{enumerate}
\end{thm} 

Perelman \cite{Pe} has announced a proof of Thurston's geometrization 
conjecture, which in turn implies the Poincar\'{e} conjecture. This 
replaces ``homotopy 3-sphere'' in the statement of the theorem by 
``3-sphere''. However, since Lopes and Morales were motivated by 
possible applications to the geometrization conjecture (see section 3) 
we do not assume geometrization in this paper.  

Lopes and Morales do not explicitly require $K$ to lie in $\Inte M$. 
The following result addresses the case where $K$ is allowed to meet 
$\bd M$. 

\begin{thm} Let $M$ be a compact, connected, orientable 3-manifold 
with $\bd M\neq\emptyset$. . 
The following are equivalent.  
\begin{enumerate}  
\item If $M_i$ is a prime connected summand of $M$, 
then $M_i$ either 
\begin{enumerate} 
\item is homeomorphic to $S^1 \times S^2$, 
\item is a homotopy 3-sphere, or 
\item has non-empty boundary. 
\end{enumerate} 
\item $\Inte M$ contains a compact, connected, irreducible 
3-manifold $V$ with $\bd V\neq\emptyset$ 
such that $\pi_1(V) \rightarrow \pi_1(M)$ is an 
isomorphism, each component of $M-\Inte V$ meets $\bd M$, and 
there is a component $W$ of $M-\Inte V$ such that any 2-sphere 
component of $\bd M$ lies in $W$, and any other component of 
$M-\Inte V$ is a collar on a component of $\bd M$. 
\item There is a knot $K$ in $M$ with $K \cap \bd M$ a non-empty 
union of disjoint arcs such that $M-K$ is irreducible 
and $\pi_1(M)$ embeds in $\pi_1(M-K)$. 
\end{enumerate}
\end{thm}

The paper is organized as follows. 
We prove Theorems 1 and 2 in sections 1 and 2, respectively. 
In section 3 we pose some questions and make some further comments.   

We work throughout in the piecewise linear category. We refer to 
\cite{He} and \cite{Ja} for basic 3-manifold topology. It will, however,  
be convenient to diverge from standard terminology by regarding every 
2-sphere embedded in a 3-manifold as being incompressible.  

\section*{\normalsize 1. Proof of Theorem 1}
\stepcounter{section}

\begin{proof}

(1)$\Rightarrow$(2): Since no $M_i$ is a 3-ball $\bd M$ contains 
no 2-spheres. 

We first show that each $M_i$ has property (2). 
In case (a) let $V=S^1\times D^2$ where $D^2$ is a closed disk in $S^2$. 
In case (b) let $V$ be a 3-ball in $M_i$. In case (c) let 
$V$ be the closure of the complement of a collar on $\bd M_i$. 

It now suffices to show that if compact, connected, 
orientable 3-manifolds $N_1$ and $N_2$ have property (2), 
then so does $N=N_1\#N_2$. If some component of $N_j-\Inte V_j$ 
is not a collar on a component of $\bd N_j$ call that component 
$W_j$. If there is no such component let $W_j$ be any component 
of $N_j-\Inte V_j$. Let $F_j=V_j\cap W_j$. Since 
$\pi_1(V)\rightarrow\pi_1(M)$ is onto $F_j$ must be connected. 
Choose a 3-ball $B_j$ 
in $\Inte N_j$ such that $B_j\cap F_j$ is a properly embedded 
disk in $B_j$ which splits $B_j$ into a 3-ball in $V_j$ 
which meets $\bd V_j$ in $B_j\cap F_j$ and a 3-ball in $W_j$ 
which meets $\bd W_j$ in $B_j\cap F_j$. Construct $N_1\#N_2$ by 
removing the interior of each $B_j$ and then gluing 
$\bd B_1$ to $\bd B_2$ so that $V_1\cap \bd B_1$ is identified with 
$V_2 \cap \bd B_2$. This creates a 3-manifold $V$ in $\Inte N$ 
which is homeomorphic to the union of $V_1$ and $V_2$ 
along disks in $F_1$ and $F_2$ and is thus irreducible. 
It also creates a component 
$W$ of $N-\Inte V$ which is homeomorphic to the 
union of $W_1$ and $W_2$ along disks in $F_1$ and $F_2$. 
and thus is not a 
collar on a component of $\bd N$ and is the only 
such component of $N-\Inte V$. Clearly $\pi_1(V)\rightarrow\pi_1(W)$ 
is an isomorphism. 

\bigskip

(2)$\Rightarrow$(3): If some component of $M-\Inte V$ is not a collar 
on a component of $\bd M$ call that component $W$. If 
there is no such component let $W$ be any component 
of $M-\Inte V$. Let $F=V\cap W$. We define submanifolds 
$X$ and $Y$ of $M$ as follows. 

If $F$ is a 2-sphere, then $V$ is a 3-ball, and we let $X=\emptyset$ and 
$Y=M$. 

If $F$ is incompressible in $V$, then we let $X=V$ and $Y=W$. 

Suppose $F$ is compressible in $V$. The compression 
creates 3-manifolds $V_1$ by cutting a 1-handle 
from $V$ and $W_1$ by adding a 2-handle to $W$. 
Note that $W_1$ is connected but might not be irreducible, 
while $V_1$ is irreducible but might not be connected. 
Let $F_1=V_1\cap W_1$. If $F_1$ is compressible in $V_1$ we 
repeat the process using a disjoint compressing 
disk to obtain $V_2$, $W_2$, and $F_2=V_2\cap W_2$. 
This procedure must stop with some $V_n$, $W_n$, and $F_n=V_n\cap W_n$. 
Each component of $F_n$ is either a 2-sphere or is 
incompressible in $V_n$. Each 2-sphere component bounds 
a 3-ball component of $V_n$. We delete these 3-balls 
from $V_n$ to obtain $X$ and add them to $W_n$ to 
obtain $Y$. In particular if $F_n$ consists entirely of 2-spheres 
then $X=\emptyset$ and $Y=M$. In any case $Y$ is connected and 
$\bd Y$ contains no 2-spheres.  

By \cite[Lemma C]{Rw} (see also \cite[Theorem 6.1]{My simple} or 
\cite[Theorem 1.1]{My excel}) there is a knot $K$ in $\Inte Y$ such 
that $Y-K$ is irreducible and $\bd Y$ is incompressible 
in $Y-K$. If $X=\emptyset$ then $M=Y$ and so $M-K$ is irreducible. 
If $X\neq\emptyset$, then $X$ is irreducible and 
$X\cap Y$ is incompressible in $X$ and $Y-K$, hence 
$M-K=X\cup(Y-K)$ is irreducible. 

By general position we can isotop $K$ in $\Inte Y$ 
so that it misses the 1-handles and 3-balls which 
were removed from $V$ in our construction. Thus $K$ 
lies in $W$, and so $V$ lies in $M-K$. Since the composition 
$\pi_1(V)\rightarrow\pi_1(M-K)\rightarrow\pi_1(M)$ is an 
isomorphism it follows 
that $\pi_1(M)$ embeds in $\pi_1(M-K)$. 

\bigskip

(3)$\Rightarrow$(1): Let $M_i$ be a prime connected summand of $M$. 
Since $M-K$ is irreducible $\bd M$ contains no 2-spheres, and so 
$M_i$ is not a 3-ball. 

We have an embedding $\pi_1(M_i)\rightarrow\pi_1(M-K)$. 
Since $\pi_1(M-K)$ has no torsion \cite[Corollary 9.9]{He} 
neither does $\pi_1(M_i)$. Hence if $\pi_1(M_i)$ is finite it 
must be trivial, and so $M_i$ is a homotopy 3-sphere 
\cite[Theorem 3.6]{He}. If $M_i$ is not irreducible, then 
it must be homeomorphic to $S^1\times S^2$ \cite[Lemma 3.13]{He}. 

Thus we may assume that $\pi_1(M_i)$ is infinite and 
$M_i$ is irreducible. It follows by the usual sphere theorem 
and Hurewicz theorem argument that $M_i$ is aspherical, 
and so $H_3(\pi_1(M_i))\cong H_3(M_i)$. Let $p:\widetilde{M-K}
\rightarrow M-K$ be the covering map with the image of the 
induced homomorphism $p_*$  
equal to the image of $\pi_1(M_i)$ in $\pi_1(M-K)$. Since 
$M-K$ is aspherical so is $\widetilde{M-K}$, and $H_3(\pi_1(M_i))
\cong H_3(\widetilde{M-K})$. Since $\widetilde{M-K}$ is non-compact 
this group must be trivial, and so $\bd M_i\neq\emptyset$. $\square$ 
\end{proof} 

\section*{\normalsize 2. Proof of Theorem 2}
\stepcounter{section}

\begin{proof} The proof follows the outline of that of Theorem 1 
with the following addtions as indicated below. 

\bigskip
 
(1)$\Rightarrow$(2): In case (c) if $M_i$ is a 3-ball, then $V_i$ 
is a 3-ball in $\Inte M_i$. 

It is easily checked that the construction of $V$ given in the 
proof of Theorem 1 has all the properties in part (2) of Theorem 
2 when $\bd M\neq\emptyset$. 

\bigskip

(2)$\Rightarrow$(3): The construction of $X$ and $Y$ is the same. 
Any 2-spheres of $\bd M$ lie in $Y$. 

Suppose the components of $Y\cap \bd M$ are $S_1,\ldots,S_p$. 
Let $\beta_0$ and $\beta_1$ be disjoint arcs in $S_1$. If $p>1$ let 
$\beta_k$ be an arc in $S_k$ for $2\leq k\leq p$. Denote the 
endpoints of $\beta_k$ by $x_k$ and $y_k$ for $0\leq k\leq p$. 
If $p=1$ let $\alpha_1$ be a properly embedded arc in $Y$ 
which joins $x_1$ to $y_0$. If $p>1$ let 
$\alpha_1,\ldots,\alpha_p$ be disjoint properly embedded arcs in $Y$ 
such that $\alpha_k$ joins $x_k$ to $y_{k+1}$ for $1\leq k\leq p-1$ and 
$\alpha_p$ joins $x_p$ to $y_0$. Let $Z$ be the exterior (closure of the 
complement of a regular neighborhood) of $\alpha_1\cup\cdots\cup\alpha_p$ 
in $Y$. Then $\bd Z$ contains no 2-spheres. Let $\alpha_0$ be a properly 
embedded arc in $Z$ which joins $x_0$ and $y_1$. 
By \cite[Proposition 6.1]{My homology}or \cite[Theorem 1.1]{My excel} 
there is a properly embedded arc $\alpha_0\p$ in $Z$ with 
$\bd\alpha_0\p=\bd\alpha_0$ such that the exterior of $\alpha_0\p$ 
in $Z$ is irreducible and has incompressible boundary. 
Let $\alpha=\alpha_0\p\cup\alpha_1\cup\cdots\cup\alpha_p$, 
$\beta=\beta_0\cup\beta_1\cup\cdots\cup\beta_p$, and $K=\alpha\cup\beta$. 
It follows as in the proof of Theorem 1 that $M-\alpha$ is irreducible 
and hence so is $M-K$. Similarly we isotop $\alpha$ into $W$ to show 
that $\pi_1(M)$ embeds in $\pi_1(M-\alpha)\cong\pi_1(M-K)$.  

\bigskip

(3)$\Rightarrow$(1): $M_i$ is now allowed to be a 3-ball. 
If it is not a 3-ball, then the proof proceeds as 
before, now using the fact that $M-\alpha$ is 
aspherical, where $\alpha$ is the closure of $K\cap\Inte M$ in $M$. $\square$ 
\end{proof}

\section*{\normalsize 3. Remarks and Questions}

We first make some observations about the Lopes-Morales examples. 
We refer to Birman \cite{Bi} for terminology. 

Let $\beta_n=\sigma_1^2\sigma_2^2\cdots\sigma_{n-1}^2$ in the 
Artin braid group $B_n$. This is a pure braid, and so its closure 
$\widehat{\beta}_n$ is a link in $S^3$ with $n$ components $K_1,\ldots,K_n$. 
Each $K_i$ is unknotted. For $i\neq j$ we have that $K_i$ and $K_j$ are 
linked if and only if $|i-j|=1$, in which case $K_i\cup K_j$ is a copy of 
the Hopf link. Let $E_n$ be the exterior of $\widehat{\beta}_n$ in 
$S^3$. Then $\bd E_n$ consists of tori $T_1,\ldots,T_n$, where 
$T_i$ is the boundary of a regular neighborhood of $K_i$ in $S^3$. 
The fact that consecutive $K_i$ are linked implies that $E_n$ is 
irreducible and $\bd E_n$ is incompressible in $E_n$.  

Lopes and Morales define $M_0$ to be $S^3$ and for $k\geq 1$ define 
$M_k$ to be $E_{2k}$. 
(Actually they use the complement $S^3-\widehat{\beta}_{2k}$ 
rather than the exterior, but we use the exterior here to make 
$M_k$ compact.) For $k\geq1$ 
they define $L_k$ to be $K_{2k-1}\cup K_{2k}$, considered as a 
link in $M_{k-1}$. They then prove algebraically that $\pi_1(M_k)$ 
embeds in $\pi_1(M_k-L_{k+1})$ for all $k\geq0$ \cite[Theorem 2.3]{LM}.

Suppose $1\leq m\leq n-1$. Let 
$\gamma_{m,n}=K_{m+1}\cup\cdots\cup K_n$. 
A closer inspection of $\widehat{\beta}_n$ 
reveals the fact that $\gamma_{m,n}$ lies in 
the interior of a collar $W$ on $\bd E_m$ in $E_m$. Let $V$ be the 
closure of $E_m-W$ in $E_m$. Since $V$ is homeomorphic to $E_m$ and 
$\bd V$ is incompressible in $W$ and hence in 
$W-\gamma_{m,n}$ we have an embedding of 
$\pi_1(E_m)$ in $\pi_1(E_m-\gamma_{m,n})$. 

Let $m=2k$. Setting $n=2k+2$ gives a new proof the Lopes-Morales result, 
and setting $n=2k+1$ gives a knot with the desired properties. 

This observation was the genesis of the construction used in the proof of 
(1)$\Rightarrow$(2)$\Rightarrow$(3) in Theorem 1.   

\bigskip

We remark that the proof of (3)$\Rightarrow$(1) in Theorems 1 and 2 
works in more general situations. In particular one does not get a 
larger class of 3-manifolds $M$ by allowing links instead of knots. 
Also, the class of $M$ in Theorem 1 is enlarged only by the inclusion of 
the class of $M$ in Theorem 2 (which is itself not enlarged) if one 
allows the embedding of $\pi_1(M)$ in a knot or link group of a different 
compact 3-manifold $N$. 

\bigskip 

Lopes and Morales were motivated by the well-known 
``$\mathbb{Q}$ conjecture'' that 
every subgroup of the fundamental group of a compact 3-manifold 
which embeds in the additive group of rational numbers must be cyclic. 
(See \cite[p. 95]{EJ}, \cite[Conjecture 1.1]{LM} and 
\cite[Conjecture 4]{My regular}.)  
This conjecture is implied by the geometrization conjecture. 
(See the references in \cite{My regular}.) On the other hand 
if it could be proven independently of the geometrization 
conjecture it might help in giving a new proof of 
the hyperbolic part of the geometrization conjecture. 
(See \cite[Corollary 2]{My regular}. It should be pointed out that 
the $\mathbb{Q}$ conjecture by itself is not known to imply the 
hyperbolization conjecture and that \cite[Conjecture 1.2]{LM} 
should have the additional hypothesis that $M$ is $P^2$-irreducible.) 

The $\mathbb{Q}$ conjecture has long been known to be true for 
Haken manifolds. (This follows easily from \cite[Corollary 3.3]{EJ}.) 
So the Lopes-Morales conjecture 
\cite[Conjecture 1.3]{LM} would imply the $\mathbb{Q}$ conjecture. 
Unfortunately Theorems 1 and 2 show that the class of orientable 
3-manifolds for which their conjecture is true does not include 
any for which the $\mathbb{Q}$ conjecture was not already known.   

\bigskip

We conclude with some questions. $M$ denotes a compact, connected 3-manifold, 
and $K$ denotes a knot in $M$.  

\begin{que} Which non-orientable $M$ contain 
a $K$ such that $\pi_1(M)$ embeds in $\pi_1(M-K)$ and $M-K$ 
is irreducible (or $P^2$-irreducible)?\end{que}

\begin{que} If $M$ contains a $K$ such that $\pi_1(M)$ embeds in 
$\pi_1(M-K)$ and $M-K$ is irreducible can one characterize all such 
knots?\end{que}

\begin{que} Which $M$ contain a $K$ such that $\pi_1(M)$ embeds in 
$\pi_1(M-K)$ and $M-K$ is hyperbolic?\end{que}

\end{document}